\documentclass[preprint,12pt]{elsarticle}

\usepackage{amssymb}
\usepackage{amsthm}
\usepackage{natbib}
\usepackage{url}
\usepackage{amsmath}
\usepackage{latexsym}
\usepackage{amsbsy}
\usepackage{slashbox}
\usepackage{caption2}
\usepackage{float}
\usepackage[all]{xy}
\usepackage[colorlinks=true, linkcolor=blue, filecolor=brown, citecolor=green,hypertex]{hyperref}
\newcommand{\set}[1]{\left\{#1\right\}}
\newcommand{\Set}[1]{\Big\{#1\Big\}}

\newcommand{\To}{\longrightarrow}
\newcommand{\X}{\mathfrak{X}}

\newcommand{\CC}{\mathbb{C}}
\newcommand{\A}{\mathcal{A}}
\newcommand{\I}{\mathcal{I}}

\newcommand{\R}{\mathcal{R}}
\newcommand{\B}{\mathcal{B}}
\newcommand{\C}{\mathcal{C}}
\newcommand{\D}{\mathcal{D}}
\newcommand{\T}{\mathcal{T}}
\newcommand{\CVR}{\mathcal{C}=(V,\mathcal{R})}
\newcommand{\p}{\mathcal{P}}
\newcommand{\x}{\chi}
\def\VRT#1{*=<5mm>[o][F-]{#1}}
     \DeclareMathOperator{\GF}{GF}
     \DeclareMathOperator{\Mat}{Mat}
   
\DeclareMathOperator{\Fib}{Fib}   \DeclareMathOperator{\Supp}{Supp}

\theoremstyle{plain}
\newtheorem{theorem}{Theorem}[section]
\newtheorem{corollary}[theorem]{Corollary}
\newtheorem{lemma}[theorem]{Lemma}
\newtheorem{proposition}[theorem]{Proposition}
\newtheorem{remark}[theorem]{Remark}
\newtheorem{prf}{Proof}
\theoremstyle{definition}
\newtheorem{definition}{Definition}[section]
\newtheorem{exmpl}[definition]{Example}

\def\thm{\begin{theorem}}
\def\ethm{\end{theorem}}
\def\prop{\begin{proposition}}
\def\eprop{\end{proposition}}
\def\rem{\begin{remark}}
\def\erem{\end{remark}}
\def\defi{\begin{definition}}
\def\edefi{\end{definition}}
\def\nmu{\begin{enumerate}}
\def\enmu{\end{enumerate}}
\def\qtn{\begin{equation}}
\def\eqtn{\end{equation}}
\def\qtnl{\begin{equation*}}
\def\eqtnl{\end{equation*}}
\def\lem{\begin{lemma}}
\def\elem{\end{lemma}}
\def\cor{\begin{corollary}}
\def\ecor{\end{corollary}}
\def\prf{\begin{proof}}
\def\eprf{\end{proof}}
\def\css{\begin{cases}}
\def\ecss{\end{cases}}
\def\exm{\begin{exmpl}}
\def\eexm{\end{exmpl}}
\journal{
Journal of Algebra.}

\biboptions{sort&compress}

\begin{document}

\begin{frontmatter}



\title{Characterization of Balanced Coherent Configurations}

\author{Mitsugu Hirasaka\corref{cor}\fnref{fn1}}
\ead{hirasaka@pusan.ac.kr}
\author{Reza Sharafdini\fnref{fn1,fn2}}
\ead{sharafdini@pusan.ac.kr}
\cortext[cor]{Corresponding Author}
\fntext[fn1]{Supported by the Korean Research Foundation grant (KRF) No: 2009-0070854.}
\fntext[fn2]{Supported by BK21 Dynamic Math Center, Department of Mathematics at Pusan National University.}
\address{Department of Mathematics, College of Natural Sciences, Pusan National University,\\ Busan 609-735, South Korea.}

\begin{abstract}
Let $G$ be a group acting on a finite set $\Omega$. Then $G$ acts on $\Omega\times \Omega$ by its entry-wise action and its orbits form the basis
relations of a coherent configuration (or shortly scheme). Our concern is to consider what follows from the assumption that the number of orbits of $G$ on $\Omega_i\times \Omega_j$ is constant whenever $\Omega_i$ and $\Omega_j$ are orbits of $G$ on $\Omega$. One can conclude from the assumption that the actions of $G$ on ${\Omega_i}$'s have the same permutation character and are not necessarily equivalent. From this viewpoint one may ask how many inequivalent actions of a given group with the same permutation character there exist. In this article we will approach to this question by a purely combinatorial method in terms of schemes and investigate the following topics: ~(i) balanced schemes and their central primitive idempotents, (ii) characterization of reduced balanced schemes.
\end{abstract}

\begin{keyword}
Coherent configurations \sep $(m,n,r)$-schemes \sep Balanced schemes \sep Central primitive idempotents.




\end{keyword}

\end{frontmatter}

\section{Introduction}

Let $G$ be a group acting on a finite set $\Omega$ with its orbits $\Omega_1,\dots,\Omega_n$ and its permutation character $\pi=\sum_{i=1}^{n}\pi_i$ where $\pi_i(g):=|\set{\alpha\in \Omega_i \mid \alpha^g=\alpha}|$ for $g\in G$. One may think what happens if $\pi_i=\pi_j$ for all $1 \leq i,j\leq n$ and can say that the number of orbits of $G$ on $\Omega_i\times \Omega_j$ by its entry-wise action is constant for all $1\leq i,j\leq n$, which motivate us to define the following concepts whose terminology is due to \cite{Inp2009}.

\defi \rm{Let $V$ be a finite set and $\R$ a set of nonempty binary relations on $V$. The pair $\CVR$ is called a
\textit{coherent configuration} (for short \textit{scheme}) on $V$ if the following conditions hold:}

\nmu[(C1)]
\item $\R$ forms a partition of the set $V\times V$.
\item $\Delta_V:=\set{(v,v)~|~ v\in V}$ is a union of  certain relations from $\R$.
\item For every $R\in\R$,\quad$R^t:=\set{(v,u)\mid(u,v)\in R}\in \R$.
\item For every $R,S,T\in\R$, the size of $\set{w\in V\mid (u,w)\in R,\ (w,v)\in S}$ does not depend on the
      choice of $(u,v)\in T$ and is denoted by $c_{RS}^{T}$.
\enmu
We say that the elements of $V$ are {\it points} and those of $\R$ are \it{basis relations}.
\edefi

\noindent Let $\CVR$ be a scheme and $\emptyset\neq X\subseteq V$. We say that $X$ is a {\it fiber} of $\C$ if $\Delta_X=\set{(x,x)\mid x\in X}\in\R$. We denote by $\Fib(\C)$ the set of all fibers of $\C$.

\defi Let $m,n$ and $r$ be positive integers. We say that a scheme $\C$ is an {\it$(m,n,r)$-scheme} if the following conditions hold:
\nmu[(i)]
\item $|\set{R\in\R\mid R\subseteq X\times Y}|=r$ for all $X,Y\in\Fib(\C)$.
\item $|X|=m$ for all $X\in\Fib(\C)$.
\item $|\Fib(\C)|=n$.
\enmu
A scheme $\C$ is  called {\it $r$-balanced} ~if {\rm (i)} holds, and {\it balanced} if it is $r$-balanced for some $r$. In Section
\ref{Section:characterization balanced} we will show that (i) implies (ii).
\edefi

Let us return to the topic in the first paragraph. Note that the orbits of $G$ on $\Omega\times \Omega$ form the basis relations of a scheme called
the \textit{2-orbit scheme} of $G$ on $\Omega$ and its fibers are $\Omega_1,\dots,\Omega_n$. Furthermore, if $\pi_i=\pi_j$ for all $1\leq i,j\leq n$,
then the  2-orbit scheme of $G$ on $\Omega$ is balanced.

We denote by  $\p(\C)$ the set of all central primitive idempotents of the adjacency algebra of $\C$ (see Section \ref{Section:Preliminaries} for details). The following theorem shows a characterization of balanced schemes in terms of their central primitive idempotents.

\thm\label{THM:main1}
Let $\C$ be a scheme. Then $\C$ is balanced if and only if for each $X\in\Fib(\C)$ the mapping  $\p(\C)\To\p(\C_X)$ {\rm($P\mapsto P_X$)} is bijective with $n_P=|\Fib(\C)|n_{P_X}$.
\ethm

One may conclude that $|\p(\C)|=r$ if $\C$ is $r$-balanced and $r\leq 5$ (see Corollary \ref{cor:r-balanced,p(C)=r}).
The following theorem deals with the converse argument for $r=1,2$.

\thm\label{THM:main2} Let $\CVR$ be a scheme. Then the following hold:
\nmu[(i)]
  \item $|\p(\C)|=1$ if and only if $\C$ is $1$-balanced.
  \item $|\p(\C)|=2$ if and only if $\C=\C_1\boxplus \C_2$ where $\C_i$ is $i$-balanced.
\enmu
\ethm

We have the following constructions of balanced schemes (see Sections  \ref{Section:characterization balanced}, \ref{Section:reduced} for the details):

\nmu[(i)]
\item Let $U$ be a union of fibers of $\C$. Then the restriction of $\C$ to $U$ is $r$-balanced if $\C$ is $r$-balanced.
\item If~$\C_i$ $(i=1,2)$ is an $(m_i,n_i,r_i)$-scheme, then $\C_1\bigotimes \C_2$ is an $\big(m_1m_2,n_1n_2,r_1r_2\big)$-scheme.
\enmu

We say that a balanced scheme $\C$ is {\it reduced} if there exist no $X,Y\in\Fib(\C)$ such that $\C_{X\cup Y}\simeq \C_X\bigotimes \T_2$ where $\T_2$
is a $(1,2,1)$-scheme (in Section \ref{Section:reduced} you will see another equivalent condition for a scheme to be reduced).
Any $r$-balanced scheme is obtained by the restriction of the tensor product of a reduced $r$-balanced scheme and a $1$-balanced scheme
(see Theorem \ref{THM:equivalence}). Now we focus our attention on reduced balanced schemes. It seems a quite difficult problem to find possible
$n$ such that there exists a reduced $(m,n,r)$-scheme for given $m$ and $r$. Actually, D.~G. Higman asked if there exists a reduced $(m,3,3)$-scheme for some $m$ (see \cite[Section 8, p.229]{Higman1987}).
Furthermore, H. Weilandt conjectured that a permutation group of prime degree has at most two inequivalent permutation representations (see \cite{Cam1972}). The following theorem shows that $n=1$ is a unique case under certain assumptions.

\thm\label{THM:main3} Let $\C$ be a reduced $(m,n,r)$-scheme and $p$ a prime. Then we have the following:
\nmu[(i)]
\item \label{THM:m<2r} If $m<2r$, then $n=1$.
\item \label{THM:p-val}  If $p\nmid m$ and $\C_X$ is $p$-valanced for some $X\in\Fib(\C)$, then $n=1$.
\enmu
\ethm

The preceding theorem is applied to characterize $(m,n,r)$-schemes up to $m\leq 11$ as follows.

\thm\label{THM:main4} If $\C$ is a reduced $(m,n,r)$-scheme and $m\leq11$, then $n\leq2$.
\ethm

Let us show the organization of this article. In Section \ref{Section:Preliminaries} we prepare some terminologies related to schemes. Section \ref{Section:characterization balanced} is devoted to balanced schemes. First we investigate the features of balanced
schemes. Indeed, we shall characterize a balanced scheme in terms of its central primitive idempotents and we prove Theorem \ref{THM:main1}.
Secondly we shall characterize schemes with at most two central primitive idempotents and we prove Theorem \ref{THM:main2}.
In Section \ref{Section:reduced} we shall extend the notion of inequivalent permutation representations to schemes. Namely, we shall define reduced $(m,n,r)$-schemes and then introduce some examples and known constructions of them to support our theory.
Finally in Section \ref{Section:Enumeration}, first we prove Theorem \ref{THM:main3}, secondly we shall enumerate reduced $(m,n,r)$-schemes for $m\leq11$ in order to prove Theorem \ref{THM:main4}.

\section{Preliminaries}\label{Section:Preliminaries}

According to \cite{Inp2009} we prepare some terminologies related to schemes. For the remainder of this section we assume that $\CVR$ is a scheme.
One can see that $V=\bigcup_{X\in\Fib(\C)}X$~(disjoint union) and
\qtn\label{qtn:partition of R}
    \R=\bigcup_{X,Y\in \text{Fib}(\C)} \R_{X,Y}\quad\quad(\text{disjoint union}),
\eqtn
where $\R_{X,Y}:= \{R\in\R\mid R\subseteq X\times Y\}$. We shall denote $\R_{X,X}$ by $\R_X$.

Let $X,Y\in\Fib(\C)$ and $R$ be a non-empty union of basis relations in $\R_{X,Y}$. For $(x,y)\in R$ we set $R_{out}(x)=\set{u\mid(x,u)\in R}$ and $R_{in}(y)=\set{v\mid(v,y)\in R}$. We shall denote the size of $R_{out}(x)$ and that of $R_{in}(y)$ by $d_R$ and $e_R$, respectively. It is easy to see that

\qtn\label{qtn:degree in&out}
    |X|d_R=|R|=|Y|e_R.
\eqtn

For $\D\subseteq\R$ we define $d_\D:=\sum_{R\in\D}d_R$ as well as $e_\D:=\sum_{R\in\D}e_R$. For instance $d_{\R_{X,Y}}=|Y|$ and $e_{\R_{X,Y}}=|X|$.
For all $X,Y\in\Fib(\C)$, we define the multi-set $d_{X,Y}:=\big\{d_R\mid R\in\R_{X,Y}\big\}$.

Note that $d_R=e_R$ for each $R\in\R$ if and only if $|X|=|Y|$ for all $X,Y\in\Fib(\C)$.
A scheme $\C$ is called \textit{half-homogeneous} if the latter condition holds.
If $\C$ is a half-homogeneous scheme, then $d_R$ ($=e_R$) is called the \textit{degree} or the \textit{valency} of $R$. Given a prime $p$ a
half-homogeneous scheme $\C$ is called {\it $p$-valenced} if the degree of each basis relation of $\C$ is a power of $p$.

A basis relation $R\in\R$ is called \textit{thin} if $d_R=e_R=1$ and a scheme $\C$ is called a \textit{homogeneous scheme} or (\textit{association scheme}) if $|\text{Fib}(\C)|=1$
or equivalently, if $\Delta_V\in \R$ (for more details regarding association schemes we refer to \cite{Zieschang2005}).
Given $X\in\Fib(\C)$ the pair $\C_X=(X,\R_X)$ is a homogeneous scheme called the \textit{homogeneous component} of $\C$ corresponding to $X$.

\defi For each $R\in\R$ we define a $\set{0,1}$-matrix $A_R$ whose rows and columns are simultaneously indexed by the elements of $V$ such that
the $(u,v)$-entry of $A_R$ is one if and only if $(u,v)\in R$. Then $A_R$ is called the {\it adjacency matrix} of $R$. Note that the subspace of
$\Mat_V(\CC)$ spanned by $\set{A_R\mid R\in\R}$ is a subalgebra, called the {\it adjacency algebra} of $\C$ and denoted by $\A(\C)$.
\edefi

Let $\A$ be the adjacency algebra of $\CVR$. Then the set $\set{A_R\mid R\in\R}$ which is a basis of $\A$ satisfies the following conditions:

\nmu
    \item[(C$^{'}$1)]  $J_V=\displaystyle\sum_{R\in\R}A_R \in \A$ where $J_V$ is the matrix whose entries are all one.
    \item[(C$^{'}$2)]  The identity matrix $I_V \in \A$.
    \item[(C$^{'}$3)]  $A_{R^t}=A_R^t$ for every $R\in\R$ where $A_R^t$ is the transpose of $A_{R}$.   %is closed under the transpose map,
    \item[(C$^{'}$4)]  For every $R, S \in \R$,\quad$A_RA_S=\displaystyle\sum_{T\in\R}c_{RS}^{T}A_T$.
\enmu

A scheme is called \textit{trivial} if all its fibers are singletons. We denote a trivial scheme on $n$ points by $\T_n$. Note that $\A(\T_n)\cong\Mat_n(\C)$ and it is easy to see that a scheme is trivial if and only if it is 1-balanced.

By $\Fib^*(\C)$ we mean the set of all non-empty unions of fibres of $\C$. Given $U\in \Fib^*(\C)$ we set $\R_U:=\set{R_U\mid R\in\R}$ where $R_U=R\cap (U\times U)$. Then the pair $\C_U=(U,\R_U)$ is a scheme on $U$ called the {\it restriction} of $\C$ to $U$. Note that $\C_U$ is homogeneous whenever $U\in\Fib(\C)$.\\
Given $U,U'\in\Fib^\ast(\C)$ we define $\A_{U,U'}$ to be the subspace of $\A$ spanned by the set $\set{A_R\mid R\in\R, R\subseteq U\times U'}$. Then the following hold:

\nmu[(i)]
\item   Given $U\in\Fib^\ast(\C)$ we have $\A_U=I_U\A I_U\cong \A(\C_U)$ where $I_U:=\sum_{\substack{X\in\Fib(\C),\\X\subseteq U}}A_{\Delta_X}$.
\item   For all $X,Y,Z,W\in\Fib(\C)$, $\A_{X,Y}\A_{Z,W}\subseteq\delta_{YZ}\A_{X,W}$.
\item   $\A=\bigoplus_{X,Y\in\Fib(\C)}\A_{X,Y}$.
\enmu

A basis relation $S$ of $\C$ is called \textit{symmetric} if $S^t=S$ and $\C$ is called \textit{symmetric}
if each basis relation of $\C$ is symmetric; and $\C$ is called \textit{commutative} if $c_{RS}^T =c_{SR}^T$ for all $R, S,T\in \R$.
This is equivalent to $A_RA_S=A_SA_R$ for all $R, S \in \R$. It is known that symmetric schemes are commutative and that the converse does not hold.
Furthermore, one can see that a commutative scheme is a homogeneous one.

\lem[{\rm \cite[(4.2)]{Higman1975}, \cite[Theorem 4.5.1]{Zieschang2005}}]\label{Lemma:r<5} If $\CVR$  is a homogeneous scheme and $|\R|\leq 5$, then $\C$ is commutative.
\elem

Given $R,S\in\R$ the \textit{complex product} of them is defined to be $RS=\big\{T\in\R\mid c_{RS}^T>0\big\}$ and  the \textit{relational product}
$R\circ S$ is defined as follows.

$$R\circ S:=\Big\{(u,v)\mid \exists~~ w\in V; (u,w)\in R, (w,v)\in S\Big\}.$$ Note that $R\circ S=\bigcup_{T\in RS}T$ and $d_{R\circ S}=d_{RS}$.

\lem\label{Lemma:constants} Let $\C$ be a scheme and $X,Y,Z\in\Fib(\C)$. Then for all $R\in\R_{X,Y}$, $S\in\R_{Y,Z}$ and $T\in\R_{X,Z}$ the following
hold:
\nmu[(i)]
\item\label{qtn:dRdS}   $d_Rd_S=\displaystyle\sum_{T\in\R_{X,Z}}c_{RS}^{T}d_T$.
\item\label{qtn:lcm}    $c_{RS}^Td_T=c_{TS^t}^Rd_R=c_{R^tT}^Sd_S$~and~${\rm lcm}(d_R,d_S)\mid c_{RS}^Td_T$.
\item\label{qtn:dgree-constant}    $d_R=\displaystyle\sum_{S\in \R_{Y,Z}}c_{RS}^T$, $e_R=\displaystyle\sum_{S\in \R_{Y,Z}}c_{R^tT}^S$, $c_{RS}^T\leq {\rm min}\set{d_R,e_S}$ and $R\R_{Y,Z}=\R_{X,Z}$.
\item\label{qtn:Kronecker} $d_R\delta_{SR^t}=c_{RS}^{\Delta_X}$ and $e_R\delta_{SR^t}=c_{SR}^{\Delta_Y}$ where $\delta$ denotes the Kronecker's delta.
\item\label{qtn:d_RS}                      $d_S\leq d_{RS}\leq d_Rd_S$ and $e_R\leq e_{RS}\leq e_Re_S$. %In particular, if $\C$ is half-homogeneous, then
%                                   ${\rm max}\set{d_R,d_S}\leq d_{RS}$.
\item\label{qtn:dT<dR}    If $RS=T$ and $d_T\leq d_R$, then $R=TS^t$.
\item\label{qtn:symmetric-dR=2} If $d_R=2$, then $RR^t=\{\Delta_X,S\}$ where $S$ is a symmetric basis relation of $\R_{X}$ with $d_S\leq2$.
\item\label{qtn:gcd}   $|RS|\leq {\rm gcd}(d_R,d_S)$.
\enmu

\elem

\prf The proof is done by the same procedure as \citep[Lemma 1.4.2, 1.4.3, 1.5.2, 1.5.3, 1.5.6]{Zieschang2005}.
\eprf

\lem\label{Lemma:d_L_S}  Let $S\in\R_{X,Y}$ and $L_S:=\set{R\in\R_X\mid RS=\set{S}}$. Then
$$d_{L_S}\mid {\rm gcd}\big(|X|,e_S\big).$$
\elem
\prf Let $y\in Y$ and $x\in  S_{in}(y)$. The condition $RS=\set{S}$ shows that $\bigcup_{R\in L_S}R_{in}(x)\subseteq S_{in}(y)$ and $\bigcup_{R\in L_S}R$ is an
equivalence relation on $X$. Since $y\in Y$ and $x\in S_{in}(y)$ are arbitrarily taken, all equivalence classes  have the same size $d_{L_S}$. It follows that $d_{L_S}$ divides both $d_S$ and $|X|$.
\eprf

\lem\label{Lemma:constant2} Let $X,Y\in\Fib(\C)$ with $X\neq Y$ and $R,S\in\R_{X,Y}$ with $R\neq S$. Then $T\in R^tR\cap S^tS$ for some $T\in\R_Y$ with
$T\neq\Delta_Y$ if and only if $c_{RS^t}^{T'}\geq2$ for some $T'\in\R_X$.
\elem

\prf Let us prove the necessity. By the assumption $c_{R^tR}^{T}\neq0$ and $c_{S^tS}^{T}\neq0$. Taking $(y,y')\in T$ (of course $y\neq y'$) there exist
$x,x'\in X$ such that $(x,y),(x,y')\in R$ and $(x',y),(x',y')\in S$. On the other hand, there exists $T'\in\R_X$ such that $(x,x')\in T'$. It follows that $c_{RS^t}^{T'}\geq2$ (see the Figure \ref{Fig:constant2}). Sufficiency follows from Figure \ref{Fig:constant2} since $c_{RS^t}^{T'}\geq2$ implies that $y\neq y'$.
\begin{figure}[H]
\begin{displaymath}
    \xymatrix{
            &\VRT{y}  \ar@{<-}[dl]_{\scriptscriptstyle\mathit{R}}\ar[dd] \ar@{<-}[dr]^{\scriptscriptstyle\mathit{S}} & \\
            \VRT{x} \ar[dr]_{\scriptscriptstyle\mathit{R}} \ar[rr]&  & \VRT{x'} \\
            &\VRT{y'}\ar@{<-}[ur]_{\scriptscriptstyle\mathit{S}}&}
\end{displaymath}\label{fig:6}
\vspace{-5mm}
\renewcommand{\captionlabeldelim}{.}
\caption{}\label{Fig:constant2}
\end{figure}
\eprf

Let $U,U'\in\Fib^*(\C)$ such that $U\cap U'=\emptyset$ and $V=U\cup U'$. Then we say that $\C$ is the {\it internal direct sum} of $\C_U$ and $\C_{U'}$  if  $|\R_{X,Y}|=1$ for all $X,Y\in\Fib(\C)$ with $X\subseteq U$ and $Y\subseteq U'$. In this case we shall write $\C=\C_U\boxplus\C_{U'}$.

Let $\C_i=(V_i,\R_i)$ $(i=1,2)$ be schemes.

We set 
$$
\R_1\otimes\R_2=\small\{R_1\otimes R_2\mid R_1\in\R_1,\ R_2\in\R_2\big\},
$$
where $R_1\otimes R_2=\Big\{\big((u_1,u_2),(v_1,v_2)\big)\mid (u_1,v_1)\in
R_1,\ (u_2,v_2)\in R_2\Big\}$. Then $\C=\big(V_1\times V_2,\R_1\otimes\R_2\big)$ is a scheme called the {\it tensor product} of $\C_1$ and $\C_2$ and
denoted by $\C_1\bigotimes\C_2$. One can see that $\Fib(\C)=\Fib(\C_1)\times\Fib(\C_2)$.

An \textit{isomorphism} from $\C_1$ to $\C_2$ is defined to be
a bijection $\psi:V_1\cup\R_1\To V_2\cup\R_2$ such that for all $u,v\in V_1$ and $R\in\R_1$, $(u,v)\in R$ if and only if $\big(\psi(u),\psi(v)\big)\in\psi(R)$. We say that $\C_1$ is \textit{isomorphic} to $\C_2$ and denote it by $\C_1\simeq \C_2$ if there exists an isomorphism from $\C_1$ to $\C_2$.

Let $\A$ be  the adjacency algebra of $\C$. Since $\A$ is closed under the complex conjugate transpose map, $\A$ is semisimple. By the Wedderburn theorem $\A$ is isomorphic to a direct sum of full matrix algebras over~$\CC$:

\qtn\label{qtn:Wedderburn}
\A=\bigoplus_{P\in\p(\C)}\A P\cong \bigoplus_{P\in\p(\C)}\Mat_{n_P}(\CC),
\eqtn
where $\p(\C)$ is the set of central primitive idempotents of $\A$, $n_P$ is a positive integer and $\Mat_{n_P}(\CC)$ is the full matrix algebra of complex $n_P\times n_P$ matrices.\\
A comparison of dimensions of the left- and right-hand sides of (\ref{qtn:Wedderburn}) shows that
\qtn
|\R|=\sum_{P\in\p(\C)}n_P^2.
\eqtn
It is known that $\C$ is commutative if and only if $n_P=1$ for each $P\in\p(\C)$.\\
Denote by $\CC^V$ the natural $\A$-module spanned by the elements of $V$. As $I_V=\sum_{P\in\p(\C)}P$ we have\\
\qtn\label{decom}
\CC^V=\bigoplus_{P\in\p(\C)}P\CC^V.
\eqtn

For each $P\in\p(\C)$ we set $m_P:=\dim_\CC(P\CC^V)/n_P$. Then the decomposition~(\ref{decom}) shows that
\qtn\label{qtn:size of V}
|V|=\sum_{P\in\p(\C)}m_Pn_P.
\eqtn
The numbers $m_P$ and $n_P$ are called the {\it multiplicity} and the {\it degree} of $P$. Set $P_0=\sum_{X} J_X/|X|$ where $X$ runs over the fibers of the
scheme~$\C$ and $J_X=\sum_{R\in\R_X}A_R$. Note that $P_0$ is a central primitive idempotent of the algebra~$\A$, which is called {\it principal}.
It is known that 

\qtn\label{qtn:m_P0n_P_0}
(m_{P_0},n_{P_0})=(1,|\Fib(\C)|). 
\eqtn
If $\C$ is homogeneous, then $P_0=J_V/|V|$ and $m_{P_0}=n_{P_0}=1$.

Below for $X\in\Fib^*(\C)$ and $P\in\mathcal{P}(\C)$ put $P_X=PI_X$ and set
$$\mathcal{P}_X(\C)=\Big\{P\in \mathcal{P}(\C)\mid\ P_X\ne 0\Big\}~~\text{and}~~ \Supp(P)=\Big\{X\in\Fib(\C)\mid\ P_X\ne 0\Big\}.$$
\thm[{\rm \cite[Proposition 2.1]{Inp1999}}] \label{THM:embeding} Let
$\CVR$ be a scheme. Then the following hold:
\nmu
 \item[{\rm~(i)}]   For each $X\in\Fib^\ast(\C)$ the mapping $P\mapsto P_X$ induces a bijection between $\mathcal{P}_X(\C)$ and $\mathcal{P}(\C_X)$.
 \item[{\rm ~(ii)}] For all $P\in\mathcal{P}(\C)$ and $X\in \Supp(P)$, $n_P=\sum_{X\in \Supp(P)}n_{P_X}$ and $m_P=m_{P_X}$ .
\enmu
\ethm

\lem\label{Lemma:union-support} Let
$\CVR$ be a scheme. Then the following hold:

\nmu[(i)]
\item $\p(\C)=\p_X(\C)$ for each $X\in\Fib(\C)$ if and only if $\Supp(P)=\Fib(\C)$ for each $P\in\p(\C)$.
\item $\Supp(P)\neq\emptyset$ for each $P\in\p(\C)$, and

\qtn\label{qtn:union-support}
           \p(\C)=\displaystyle\bigcup_{X\in\Fib(\C)}\p_X(\C).
\eqtn
Besides, $\p(\C)=\p_U(\C)\cup\p_{U'}(\C)$ where $U,U'\in\Fib^*(\C)$ with $U\cap U'=\emptyset$ and $V=U\cup U'$.
\enmu

\elem

\prf (i) Let $X\in\Fib(\C)$ and $P\in\p(\C)$. Then $P\in\p_X(\C)$ if and only if $X\in\Supp(P)$.
This completes the proof.\\
(ii) Let $P\in\p(\C)$ such that $\Supp(P)=\emptyset$. Then for all $X\in\Fib(\C)$, $PI_X=0$ and then
$P=PI_V=\sum_{X\in\Fib(\C)} PI_X=0$, a contradiction. Therefore, $\Supp(P)\neq\emptyset$. Let $P\in\p(\C)$, as $\Supp(P)\neq\emptyset$, there exists
$X\in\Fib(\C)$ such that $PI_X\neq0$. This means that $P\in\p_X(\C)$ and the proof of (\ref{qtn:union-support}) is completed.

Let $P\in\p(\C)$. Then $P\in\p_X(\C)$ for some $X\in\Fib(\C)$. Since $V=U\cup U'$, $X\subseteq U$ or $X\subseteq U'$. It follows that $P\in\p_U(\C)$ 
or $P\in\p_{U'}(\C)$. This completes the proof. 
\eprf

\prop [{\rm \cite[p.223]{Higman1987}, \cite[p.22 (8.1)]{Higman1975}}]\label{Prop:Higman-d_XY} Let $\CVR$ be a scheme.
Then the following hold:

\nmu[(i)]
\item  Let $X,Y\in\Fib^*(\C)$ such that $X\cap Y=\emptyset$ and $V=X\cup Y$. Then
         $$\dim_{\CC}(\A_{X,Y})=\sum_{P\in \p_X\bigcap \p_Y}n_{P_X}n_{P_Y}.$$
\item  For all $X,Y\in\Fib(\C)$, $|\R_{X,Y}|=\displaystyle\sum_{P\in \p_X\bigcap \p_Y}n_{P_X}n_{P_Y}$.
\enmu

\eprop

\lem\label{Lemma:direct sum1}
Let $\CVR$ be a scheme with the adjacency algebra $\A(\C)$. If $U,U'\in\Fib^\ast(\C)$ such that $U\cap U'=\emptyset$, then

\qtn\label{qtn:direct sum1}
|\Fib(\C_U)||\Fib(\C_{U'})|\leq \dim_{\CC}(\A_{U,U'}).
\eqtn
Furthermore, the equality holds if and only if $\C_{U\cup U'}=\C_U\boxplus\C_{U'}$.
\elem

\prf It is clear that $\Fib(\C_U)=\set{X\in\Fib(\C)\mid X\subseteq U}$ as well as $\Fib(\C_{U'})=\set{Y\in\Fib(\C)\mid Y\subseteq U'}$. By the assumption,

$$
\A_{U,U'}=\bigoplus_{\substack{X\in\Fib(\C_U),\\Y\in\Fib(\C_{U'})}}\A_{X,Y}.
$$

On the other hand, $\dim_{\CC}(\A_{X,Y})=|\R_{X,Y}|\geq1$ for each $X\in\Fib(\C_U)$ and $Y\in\Fib(\C_{U'})$. It follows that

$$|\Fib(\C_U)||\Fib(\C_{U'})|\leq\sum_{\substack{X\in\Fib(\C_U),\\Y\in\Fib(\C_{U'})}}|\R_{X,Y}|=\dim_{\CC}(\A_{U,U'}).$$

The equality holds if and only if $|\R_{X,Y}|=1$ for all $X\in\Fib(\C_U)$ and $Y\in\Fib(\C_{U'})$ which is exactly the definition of internal direct sums.

\eprf

\lem\label{Lemma:direct sum2} Let $\CVR$ be a scheme with the principal idempotent $P_0$ and let $U,U'\in\Fib^\ast(\C)$ such that $U\cap U'=\emptyset$ and $V=U\cup U'$.
Then $\C=\C_U\boxplus\C_{U'}$ if and only if
$\p_U(\C)\cap\p_{U'}(\C)=\set{P_0}$. %Furthermore, $|\p(\C)|=|\p(\C_U)|+|\p(\C_U)|-1$
\elem

\prf  Let us prove the sufficiency first. It is clear that $P_0\in \p_U(\C)\cap\p_{U'}(\C)$. By Lemma \ref{Lemma:direct sum1} and Proposition \ref{Prop:Higman-d_XY}~(i)
we have $$|\Fib(\C_U)||\Fib(\C_{U'})|=\sum_{P\in \p_U\cap\p_{U'}}n_{P_U}n_{P_{U'}}.$$

Since $n_{P_{0U}}=|\Fib(\C_U)|$ and $n_{P_{0U'}}=|\Fib(\C_{U'})|$, it follows that
\qtn
\p_U(\C)\cap\p_{U'}(\C)=\set{P_0}.
\eqtn

Conversely, if $\p_U(\C)\cap\p_{U'}(\C)=\set{P_0}$, then by Proposition \ref{Prop:Higman-d_XY}~(i),  $$\dim_{\CC}(\A_{U,U'})=|\Fib(\C_U)||\Fib(\C_{U'})|.$$ It follows from Lemma \ref{Lemma:direct sum1} that $\C=\C_U\boxplus\C_{U'}$.

\eprf

\section{Characterization of balanced schemes}\label{Section:characterization balanced}

\noindent\textbf{Proof of Theorem \ref{THM:main1}:}\\
First we prove the necessity. Let $X,Y\in\Fib(\C)$. By Proposition \ref{Prop:Higman-d_XY},
$|\R_{X,Y}|=\sum_{P\in{\p_X\cap\p_Y}} n_{P_X}n_{P_Y}$. By the Cauchy-Schwarz inequality we have

\begin{eqnarray*}
|\R_{X,Y}|^2=\big(\sum_{P\in \p_X\cap\p_Y} n_{P_X}n_{P_Y} \big)^2\label{Cauchy-Schwarz}
&\leq& \sum_{P\in \p_X\cap\p_Y} n_{P_X}^2 \sum_{P\in\p_X\cap\p_Y}n_{P_Y}^2\\
&\leq& \sum_{P\in\p_X} n_{P_X}^2 \sum_{P\in\p_Y}n_{P_Y}^2\\
&=&|\R_{X}||\R_{Y}|=|\R_{X,Y}|^2.
\end{eqnarray*}

This implies that
$$\big(\sum_{P\in{\p_X\cap\p_Y}} n_{P_X}n_{P_Y}\big)^2=\sum_{P\in\p_X} n_{P_X}^2\sum_{P\in\p_Y}n_{P_Y}^2.$$

It follows that $\p_X(\C)=\p_Y(\C)$ and thus applying Lemma \ref{Lemma:union-support}~(i) we have $\p(\C)= \p_X(\C)$. Consequently, the mapping $\p(\C)\rightarrow \p(\C_X)~(P\mapsto P_X)$ is well-defined and bijective by Theorem \ref{THM:embeding}. Since the equality holds in the Cauchy-Schwarz inequality,
we have $\Big\langle n_{P_X}|P\in\p(\C)\Big\rangle=\alpha\Big\langle n_{P_Y}|P\in\p(\C)\Big\rangle$.
However, $\alpha=1$ since $|\R_X|=|\R_Y|$. Hence, $n_{P_X}=n_{P_Y}$ for all $P\in\p(\C)$. Therefore, by Theorem \ref{THM:embeding} and
Lemma \ref{Lemma:union-support}~(ii),
$$n_P=\sum_{X\in\Supp(P)}n_{P_X}=\sum_{X\in\Fib(\C)}n_{P_X}=|\Fib(\C)|n_{P_X}.$$
Now let us prove the sufficiency. Given $X,Y\in\Fib(\C)$ the assumption along with Theorem \ref{THM:embeding} assert that $\p(\C)=\p_X(\C)=\p_Y(\C)$ and $n_{P_X}=n_{P_Y}$ for each $P\in\p(\C)$. On the other hand, by Preposition \ref{Prop:Higman-d_XY}~(ii), we have
$$|\R_{X,Y}|=\displaystyle\sum_{P\in{\p_X\cap\p_Y}} n_{P_X}n_{P_Y}=
\sum_{P\in\p(\C)} n_{P_X}^2=|\R_X|.$$
Hence, $\C$ is balanced.\hspace{9.8cm}$\square$

\cor\label{cor:r-balanced,p(C)=r} Let $\C$ be an $r$-balanced scheme.  If $\C_X$ is commutative for some $X\in\Fib(\C)$, then $|\p(\C)|=r$. For instance, $|\p(\C)|=r$ if $r\leq 5$.
\ecor

\prf Let $X\in\Fib(\C)$. Since $\C_X$ is commutative, $|\p(\C_X)|=|\R_X|=r$. By Theorem \ref{THM:main1},
$|\p(\C)|=|\p(\C_X)|=r$. In particular, if $r\leq 5$, then by Lemma \ref{Lemma:r<5}, $\C_X$ is commutative and thus $|\p(\C)|=r$.
\eprf

\noindent\textbf{Proof of Theorem \ref{THM:main2}~(i):}\\
Let $X\in\Fib(\C)$. By Theorem \ref{THM:embeding}, $|\p(\C_X)|=1$. On the other hand, $\C_X=(X,\R_X)$ is a homogeneous scheme, so $|X|=m_{P_{0X}}n_{P_{0X}}=1$, by (\ref{qtn:m_P0n_P_0}).
Hence, every fiber of $\C$ is a singleton and thus $\C$ is trivial. Conversely, the adjacency algebra of a trivial scheme is the full matrix algebra and
thus it has  only one central primitive idempotent.\hspace{+8cm}$\square$\\

In order to prove Theorem \ref{THM:main2}~(ii), we need the following preparations.

\lem[{\rm \cite[(4.2)]{Higman1975}}]\label{Lemma:p=2homoge} Let $\CVR$ be a homogenous scheme. Then $|\p(\C)|=2$ if and only if $|\R|=2$.
\elem

\lem\label{Lemma:p=2} Let $\CVR$ be a scheme. If $\p(\C)=\set{P_0,P_1}$ with $P_0\neq P_1$, then the following hold:

\nmu
      \item [\rm~(i)] $X\notin \Supp(P_1)$ if and only if $|X|=1$.
      \item [\rm~(ii)] $|\R_X|=
                                \css
                                     2 & \quad  \text{if}~~X\in \Supp(P_1) \\
                                     1 & \quad  \text{if}~~X \notin \Supp(P_1)
                                \ecss.$
\enmu

\elem

\prf ~(i) Since $I_V=P_0+P_1$, $P_1=\sum_{X\in\Fib(\C)}(I_X-J_X/|X|)$. Let $X\in\Fib(\C)$. Then $X\notin \Supp(P_1)$ if and only if $0=P_1I_X=I_X-J_X/|X|$ if and only if $|X|=1$.\\
~(ii) If $X\in \Supp(P_1)$, then $P_{1}I_X\neq0$ and  by
Theorem \ref{THM:embeding}, $|\p(\C_X)|=2$. Since $\C_X=(X,\R_X)$ is homogeneous, it follows from Lemma \ref{Lemma:p=2homoge} that $|\R_X|=2$.
If $X\notin \Supp(P_1)$, then by (i), we have $|X|=1$. It follows that
$|\R_X|=1$.
\eprf
\lem\label{Lemma:half2} Let $\CVR$ be a scheme. If $\p(\C)=\set{P_0,P_1}$ with $P_0\neq P_1$, then the following hold:

\nmu
    \item [\rm~(i)]   $n_{P_1}=|\Supp(P_1)|$ and $|X|=1+m_{P_1}$ for each $X\in \Supp(P_1)$.
    \item [\rm~(ii)]  $|\R_{X,Y}|=2$ for each $X,Y\in\Supp(P_1)$.
\enmu
\elem

\prf (i) Let $X\in\Supp(P_1)$. By Lemma \ref{Lemma:p=2}, $|\R_X|=2$ and thus by Lemma \ref{Lemma:r<5}, $\C_X$ is commutative. By Theorem \ref{THM:embeding} we have

$$n_{P_1}=\sum_{X\in\Supp(P_{1})}n_{P_{1X}}=|\Supp(P_{1})|.$$
Thus (\ref{qtn:size of V}) implies that $|X|=1+m_{P_1}$.\\
\rm~(ii) Let $X,Y\in\Supp(P_1)$. Then by Lemma \ref{Lemma:p=2}~(ii), $|\R_Y|=|\R_X|=2$. It follows from Lemma \ref{Lemma:p=2homoge} that
$\p_X(\C)\cap\p_Y(\C)=\p(\C)$. Therefore, Proposition \ref{Prop:Higman-d_XY}~(ii)
implies that $|\R_{X,Y}|=2.$
\eprf

\noindent\textbf{Proof of Theorem \ref{THM:main2}~(ii):} Let $\p(\C)=\set{P_0,P_1}$ and set $U:=\bigcup_{X\in\Supp(P_1)}X$ and $U':=V\setminus U.$
If $X\in \Supp(P_1)$ and $Y\notin \Supp(P_1)$,
then $|\R_{X,Y}|=1$, since $|Y|=1$ by Lemma \ref{Lemma:p=2}. This implies that
$$\C=\C_U\boxplus \C_{U'}.$$
Note that by Lemma \ref{Lemma:p=2}~(ii) and Lemma \ref{Lemma:half2}~(ii), $\C_U$ is 2-balanced and $\C_{U'}$ is 1-balanced.
Conversely, by Lemma \ref{Lemma:direct sum2} and  Corollary \ref{cor:r-balanced,p(C)=r},
$|\p(\C)|=|\p(\C_1\boxplus\C_2)|=|\p(\C_1)|+|\p(\C_2)|-1=|\p(\C_2)|=2$.\hspace{+4cm}$\square$\\

\cor\label{Corollary:half-homo} If $\C$ is a balanced scheme, then the following hold:

\nmu[(i)]
\item   $\C$ is half-homogeneous.
\item   For every $X,Y\in\Fib(\C)$, $\A_X$ and $\A_Y$ are isomorphic as $\CC$-algebras.
\enmu

\ecor

\prf
(i) Let $X\in\Fib(\C)$ and consider the scheme $\C_X=(X,\R_X)$. It follows from Theorem \ref{THM:main1} that the mapping $\p(\C)\To \p(\C_X)$
{\rm($P\mapsto P_X$)} is bijective with $n_P=|\Fib(\C)|n_{P_X}$. By (\ref{qtn:size of V}) and Theorem \ref{THM:embeding}~(ii), the size of $X$ is computed as follows.

\begin{eqnarray*}
|X|=\sum_{P\in\p(\C)}n_{P_X}m_{P_X}=\frac{1}{|\Fib(\C)|}\sum_{P\in\p(\C)}n_{P}m_{P}=\frac{|V|}{|\Fib(\C)|}.
\end{eqnarray*}

Hence, the size of each fiber is constant and thus $\C$ is half-homogeneous.\\

(ii) Let $X,Y\in\Fib(\C)$ and $P\in\p(\C)$. Then by Theorem \ref{THM:main1}, $n_{P_X}=n_{P_Y}$. It follows from (\ref{qtn:Wedderburn}) that

$$\A_X=\bigoplus_{P\in\p(\C)}\Mat_{n_{P_X}}(\CC)\cong\bigoplus_{P\in\p(\C)}\Mat_{n_{P_Y}}(\CC)=\A_Y.$$

\eprf

%\defi An $(m,n,r)$-scheme $\C$  is called {\it trivially decomposable} if it is the tensor product of
%an $(m,1,r)$-scheme and the $(1,n,1)$-scheme.
%\edefi

Given a scheme $\C$ we define a relation $E_{\C}$ on $\Fib(\C)$ as follows.
\qtn\label{qtn:EC}
E_{\C}:=\Set{(X,Y)\in\Fib(\C)\mid \exists R\in \R_{X,Y}~;~d_R=e_R=1}.
\eqtn
\lem $E_{\C}$ is an equivalence relation on $\Fib(\C)$.
\elem

\prf For each $X\in\Fib(\C)$, $\Delta_X$ is a thin basis relation in $\R_X$ and thus $E_{\C}$ is reflexive. If $R\in\R_{X,Y}$ is thin,
then $R^t\in\R_{Y,X}$ is also thin and then $E_{\C}$ is symmetric. Let $X,Y,Z\in\Fib(\C)$ and $R\in \R_{X,Y},  S\in \R_{Y,Z}$ such that $d_R=d_S=1$
and $e_R=e_S=1$. It follows from Lemma \ref{Lemma:constants}~(\ref{qtn:d_RS}) that $RS$ is a thin basis relation in $\R_{X,Z}$ and thus $E_{\C}$ is
transitive.
\eprf

\thm\label{THM:equivalence} Let $\C$ be an $(m,n,r)$-scheme with the equivalence relation $E_{\C}$. If $\X$ is a transversal of $E_{\C}$ in $\Fib(\C)$, then $\C$ is a restriction of $\C_{U_{\X}}\bigotimes\T_n$ where $U_{\X}=\bigcup_{X\in\X}X$ and $\T_n$ is a $(1,n,1)$-scheme.
\ethm

\prf

Let $I_n:=\set{1,\dots,n}$ and $E_n:=\set{e_{ij}~\mid~ 1\leq i,j \leq n}$ where for $1\leq i,j\leq n$, $e_{ij}=\set{(i,j)}$. Then  $\T_n=(I_n,E_n)$.
Let $\X=\set{X_1,\dots, X_s}$ and suppose that for each $i\in\set{1,2,\dots,s}$, $E_{\C}(X_i)=\set{X_{i1},X_{i2},\dots,X_{im_i}}$ where
$X_{i1}:=X_i$ and $X_{ij}$'s are distinct fibers. In this case, $V=\bigcup_{i=1}^{s}\bigcup_{j=1}^{m_i}X_{ij}$.
%Let $U_{i}:=\bigcup_{j=1}^{m_i}X_{ij}$. Then, $\C_{U_{i}}$ is trivially decomposable.
For all $i\in\set{1,\dots,s}$ and $j\in\set{1,\dots,m_i}$, there exists $R_{ij}\in\R_{X_{i1},X_{ij}}$ with $d_{R_{ij}}=1$.
Therefore, there exists a bijection $R_{ij}:X_{i}\To X_{ij}, (x_i\mapsto x)$ where $x$ is the unique element of $X_{ij}$ such that $(x_i,x)\in R_{ij}$.
Indeed, $R_{ij}(X_{i})=X_{ij}$. Thus, for each $x\in V$, there exist unique $i\in\set{1,\dots,s}$ and $j\in\set{1,\dots,m_i}$ such that
$R_{ij}(X_{i})=X_{ij}$ and $x\in R_{ij}(X_{i})$. We claim that the map $\psi$ defined as follows is a monomorphism.
\begin{eqnarray*}
\psi:V\cup\R&\To&\big(U_{\X}\times I_n\big)\cup(\R_{U_{\X}}\bigotimes E_n).\\
x&\mapsto& (x_i,j);        \quad   \quad R_{ij}(x_i)=x,   \\
R&\mapsto& R_{ij}RR_{kl}^t\otimes e_{jl}; \quad R\in\R_{X_{ij},X_{kl}}
\end{eqnarray*}

Note that $\psi$ is injective, since $R_{ij}$ is a bijection for all $i\in\set{1,\dots,s}$ and $j\in\set{1,\dots,m_i}$.
Let $(x,y)\in R$ and $R\in\R_{X_{ij},X_{kl}}$. Then there exists $(x_i,y_k)\in X_i\times X_k$ such that $R_{ij}(x_i)=x$ and $R_{kl}(y_k)=y$.
This means that $(x_i,y_k)\in R_{ij}RR_{kl}^t$. It follows that, $(\psi(x),\psi(y))=((x_i,j),(y_k,l))\in\psi(R)$.
%Conversely, if $(x_i,y_k)\in R_{ij}RR_{kl}^t$, then there exists $(u,v)\in R\in\R_{X_{ij},X_{kl}}$ such that $(x_i,u)\in R_{ij}$ and $(y_k,v)\in R_{kl}$.
%However, $(u,v)=(x,y)$, since $R_{ij}$ is a bijection. Hence, $\psi$ is a monomorphism.

\eprf
The following is an immediate consequence of the preceding theorem.

\cor\label{Cor:main T.D}
 Let $\CVR$ be an $(m,n,r)$-scheme. Then $\C\simeq \C_X\otimes \T_n$ for $X\in\Fib(\C)$ if and only if $E_{\C}$ is trivial.
\ecor

\section{Reduced $(m,n,r)$-schemes}\label{Section:reduced}

%In this section we shall extend the notion of inequivalent permutation representations to schemes. To do so, we shall define reduced $(m,n,r)$-schemes. Later we shall cite some examples and known constructions of reduced $(m,n,r)$-schemes to support our theory.

\defi An $(m,n,r)$-scheme $\C$ is called \textit{reduced }if its equivalence relation $E_{\C}$ is discrete.
\edefi
\rem Note that by Corollary \ref{Cor:main T.D}, a balanced scheme $\C$ is {\it reduced} if and only if there exist no $X,Y\in\Fib(\C)$ such that $\C_{X\cup Y}\simeq \C_X\bigotimes \T_2$ where $\T_2$ is a $(1,2,1)$-scheme.
\erem
\exm\label{Example:Symmetric design} Let $(X,\B,\I)$ be a symmetric design with the set $X$ of points,  the set $\B$ of blocks and the incidence relation $\I\subseteq X\times\B$. Set $V=X\cup \B$ (disjoint union) and define the relations
$R_i~(i=1,\ldots,8)$ on $V$ as follows.

\begin{align*}
R_1&=\Delta_X,\quad & R_2&=\Delta_{\B}, &
R_3&=(X\times X)\setminus\Delta_X,\quad & R_4&=(\B\times \B)\setminus\Delta_\B,\\
R_5&=\I,\quad & R_6&=R_5^t,\quad &
R_7&=(X\times \B)\setminus\I,\quad & R_8&=R_7^t.
\end{align*}

It is known that $(V,\set{R_i}_{i=1}^{8})$ is an $(m,2,2)$-scheme where $m=|X|$. %Note that, if a symmetric design on $m$ points is non-trivial, then the corresponding $(m,2,2)$-scheme is reduced.
\eexm

%\exm The following  are two well-known examples of reduced schemes.
%\nmu[(i)]
%\item A non-trivial linked symmetric design is identified with an $(m,n,2)$-scheme for some $m,n$  {\rm(see \cite{Cam1972})}.
%\item An strongly regular design of the second kind corresponds to a reduced $(m,2,3)$-scheme {\rm(see \cite{Higman-SRD95})}.
%\enmu
%\eexm

In \cite{Cam1972} and \cite{Higman-SRD95} they mentioned $(m,n,2)$-schemes and $(m,2,3)$-schemes as linked symmetric designs and
strongly regular designs of the second kind, respectively.

%Recall that a permutation representation of a group $G$ on $V$ is called {\it doubly transitive} if the number of orbits of $G$ on $V\times V$, under its
%entry-wise action, is exactly 2. %The following example gives a construction of linked symmetric designs.

\exm[{\rm \cite{Polla}, \cite[p.6, Example~(i)]{Cam1972}}]
 The split extension of the translation group of the vector space $\GF(2^t)^{2k}$ by the
symplectic group ${\rm Sp}(2k,2^t)$ has $2^t$ pairwise inequivalent doubly transitive representations of degree $2^{2kt}$ with the same character.
This gives a reduced $(2^{2kt},2^t,2)$-scheme.
\eexm

\lem[{\rm \cite[Chapter 4]{Cam2001}}]\label{Lemma:Cameron}
Let $G={\rm PGL}(t,q)$ and $\Omega_k$ the set of $k$-dimensional subspaces of  the vector space ${\GF}(q)^t$.
Let $\pi_k$ denote the permutation character of $G$ on $\Omega_k$.  Then we have the following:\\

For each $k\leq \frac{t}{2}$ there exist irreducible characters $\x_0,\x_1,\dots,\x_k$  of $G$ with $\x_0 =1_G$ such that

\qtn\label{k=r-1}
\pi_{t-k}=\pi_k=\sum_{i=0}^{k}\x_i.
\eqtn
\elem

\exm\label{Cor:existance}
Let $r$ and $t$ be positive integers such that $r-1\leq \frac{t}{2}$. Then by Lemma \ref{Lemma:Cameron}, the 2-orbit scheme of~${\rm PGL}(t,q)$ on $\Omega_{r-1}\cup \Omega_{t-r+1}$ is a reduced $\big(\genfrac{[}{]}{0pt}{}{t}{r-1}_{q},2,r\big)$-scheme, say $\C$. Moreover, as the actions of ${\rm PGL}(t,q)$ on both $\Omega_{r-1}$ and $\Omega_{t-r+1}$ are multiplicity free, both $\C_{\Omega_{r-1}}$ and $\C_{\Omega_{t-r+1}}$ are commutative and hence by Corollary \ref{cor:r-balanced,p(C)=r}, $|\p(\C)|=r$.
\eexm

\lem\label{Lemma:reduced construct} Let $\C_i$ be an $(m_i,n_i,r_i)$-scheme for $i=1,2$.
Then $\C_1\bigotimes\C_2$ is an $\big(m_1m_2,n_1n_2,r_1r_2\big)$-scheme. Furthermore, $\C_1\bigotimes\C_2$ is reduced if
and only if both $\C_1$ and $\C_2$ are reduced.
\elem

\prf The first statement is obtained by the definition of $\C_1\bigotimes\C_2$. Let $R_i$ be a basis relation of $\C_i$ for $i=1,2$. Then $R_1\otimes R_2$ is thin if and only if  both $R_1$ and
$R_2$ are thin. This implies that $\C_1\bigotimes\C_2$ is reduced if and only if both $\C_1$ and $\C_2$ are reduced.
\eprf

\noindent\textbf{Problem 1:} Given an odd prime $p$ does there exist any reduced $(m,3,p)$-scheme for some $m$?\\

The following problem is inspired from a conjecture by H. Weilandt on permutation representations (see \cite{Cam1972}, Remark \ref{Remark:(m,n,2)}
and Lemma \ref{Lemma:(m,n,2)a,b}).\\
\noindent\textbf{Problem 2:}
If $\C$ is a reduced $(p,n,r)$-scheme for some $r$ and prime $p$, then $n\leq2$.

%If $\C$ is a reduced $(m,2,3)$-scheme, then $\C_X$ is symmetric for each $X\in\Fib(\C)$
%(see Lemma \ref{Lemma:Higman symmetric}).\\
%\noindent\textbf{problem 3:}
%If $\C$ is a reduced $(m,n,p)$-scheme for some $m$ and prime $p$, then $\C_X$ is symmetric for each $X\in\Fib(\C)$.

\section{Enumeration of $(m,n,r)$-schemes for $m\leq 11$}\label{Section:Enumeration}

\noindent\textbf{Proof of Theorem \ref{THM:main3}~(i):} Let $\C$ be a reduced $(m,n,r)$-scheme and $X,Y\in\Fib(\C)$ with $X\neq Y$. Then $2\leq d_R$ for each
$R\in\R_{X,Y}$ and
$$2|\R_{X,Y}|\leq\sum_{R\in\R_{X,Y}}d_R=m,$$
a contradiction.

In order to prove Theorem \ref{THM:main3}~(ii) we need the following lemma.
\lem\label{Lemma:coprime} Let $\C$ be an $(m,n,r)$-scheme and $X,Y,Z\in\Fib(\C)$. If $T\in\R_{X,Y}$ such that $d_T$ is prime to
$\prod_{R\in\R_{Y,Z}} d_R$, then $d_{Y,Z}$ coincides with $d_{X,Z}$ as multi-sets and $d_T\leq\min\big\{d_R\mid R\in \R_{Y,Z}\big\}$.
\elem

\prf For each $R\in\R_{Y,Z}$, ${\rm gcd}(d_T,d_R)=1$. By Lemma \ref{Lemma:constants}~(\ref{qtn:gcd}), $|TR|=1$ and we may define the following map.
\begin{eqnarray*}
\psi:\R_{Y,Z}&\To&\R_{X,Z}\\
R&\mapsto& S;\quad TR=\set{S}.
\end{eqnarray*}

By Lemma \ref{Lemma:constants}, $\psi$ is surjective. Since $|\R_{Y,Z}|=|\R_{X,Z}|$, $\psi$ must be a bijection.
Consequently, $\sum_{R\in\R_{Y,Z}}d_R=\sum_{S\in\R_{X,Z}}d_{S}=\sum_{R\in\R_{Y,Z}}d_{TR}$. On the other hand, by Lemma \ref{Lemma:constants}~(\ref{qtn:d_RS}), $
d_R\leq d_{TR}$ for each $R\in\R_X$ and thus $d_R=d_{TR}$ for each $R\in\R_{Y,Z}$. Furthermore, by Lemma \ref{Lemma:constants}~(\ref{qtn:d_RS}),
$d_T\leq d_{TR}=d_R$ for each $R\in\R_{Y,Z}$.
\eprf

\noindent\textbf{Proof of Theorem \ref{THM:main3}~(ii):} Let $\C$ be a reduced $(m,n,r)$-scheme and let $X\in\Fib(\C)$ such that $\C_X$ is $p$-valanced. Clearly $m=\sum_{T\in\R_{X,Y}}d_T$ where $X,Y\in\Fib(\C)$ with $X\neq Y$. Since $p\nmid m$, so there exists $T\in\R_{X,Y}$ such that $p\nmid d_T$. Since $\C_X$ is $p$-valenced, $d_T$ is prime to
$\prod_{R\in\R_{X}} d_R$. As $d_{\Delta_X}=1$, it follows from Lemma \ref{Lemma:coprime} that
$d_T\leq\min\big\{d_R\mid R\in \R_X\big\}=1,$
a contradiction.~~~~~\hspace{12cm}$\square$
\lem\label{Lemma:(m,n,2)symmetric} Let $\C$ be an $(m,n,2)$-scheme and $R\in\R_{X,Y}$ where $X,Y\in\Fib(\C)$. Then $d_R(d_R-1)=\lambda(m-1)$ for some $\lambda\in \mathbb{N}$.
\elem

\prf
Let $\C$ be an $(m,n,2)$-scheme and $X,Y\in\Fib(\C)$. For each $R\in\R_{X,Y}$ we have by Lemma \ref{Lemma:constants}~(\ref{qtn:dRdS}),
$$A_RA_{R^t}=\sum_{S\in\R_X}c_{RR^t}^S A_S=d_R I_X+ c_{RR^t}^{\Delta_X^c}(J_X-I_X),$$
where $\Delta_X^c=(X\times X)\setminus\Delta_X$. It follows that $R\in\R_{X,Y}$ is regarded as the incident relation of a symmetric
$(m, d_R,\lambda)$-design where $\lambda=c_{RR^t}^{\Delta_X^c}$. A basic property of symmetric deigns implies that $d_R(d_R-1)=\lambda(m-1)$.
\eprf

\rem\label{Remark:(m,n,2)} Let $m$ and $t$ be positive integers and  $q$ an odd prime power such that

\qtn\label{qtn:(m,n,2)}
m-1=2^tq.
\eqtn
Then there are exactly four $d\in\set{1,\dots,m-1}$ such that $d(d-1)\equiv0$ $({\rm mod} ~2^tq)$ by an elementary number theoretical argument.
It follows that if $\C$ is a reduced $(m,n,2)$-scheme,  then $d_{X,Y}$ is uniquely determined for all $X,Y\in\Fib(\C)$. In particular, if m is prime
satisfying (\ref{qtn:(m,n,2)}), then there is $\gamma \in \set{2,\dots,m-2}$ such that $d_{X,Y}=\set{\gamma, m-\gamma}$ and ${\rm gcd}(\gamma,m-\gamma)=1$ for all $X,Y\in\Fib(\C)$.

\erem

\lem\label{Lemma:(m,n,2)a,b} Let $\C$ be a reduced $(m,n,2)$-scheme. Suppose that $d_{X,Y}=\set{a,b}$ with ${\rm gcd}(a,b)=1$ for all $X,Y\in\Fib(\C)$.
Then $n\leq2$.
\elem

\prf Suppose that $X,Y$ and $Z$ are distinct fibers of $\C$ and let $\R_{X,Y}=\set{R,R'}$, $\R_{Y,Z}=\set{S,S'}$,
$\R_{X,Z}=\set{T,T'}$ so that $d_R=d_S=d_T=a<b=d_{R'}=d_{S'}=d_{T'}$.
By Lemma \ref{Lemma:constants}~(\ref{qtn:dRdS}),~(\ref{qtn:lcm}),~(\ref{qtn:dgree-constant}), $a^2=d_Rd_S=\alpha a+\beta b$ such that $a\mid b\beta$ and
$\beta <a$. Since ${\rm gcd}(a,b)=1$, it follows that $\beta=0$ and $\alpha=a$. This implies that $c_{RS}^{T}=a=d_R$. It follows that
%Let $x\in X$. Then $R_{out}(x)=R_{out}(x)\cap S_{in}(z)$ for all $z\in T_{out}(x)$, and hence,
\qtn\label{qtn:(a,b)=1}
R_{out}(x)\subseteq S_{in}(z),
\eqtn
where $(x,z)\in T$. Now we take  $y_1,y_2\in R_{out}(x)$ so that $y_1\neq y_2$. It follows from (\ref{qtn:(a,b)=1}) that $T_{out}(x)\subseteq S_{out}(y_1)\cap S_{out}(y_2)$.
This fact along with Lemma \ref{Lemma:constants}~(\ref{qtn:dgree-constant}) assert that $a=c_{SS^t}^{\Delta_Y^c}$ where
$\Delta_Y^c=(Y\times Y)\setminus\Delta_Y$. Therefore, by Lemma \ref{Lemma:(m,n,2)symmetric}, $a(a-1)=a(a+b-1)$. It follows that $ab=0$, a contradiction. This completes the proof.
\eprf

\lem \label{Lemma:r=2}Let $\C$ be a reduced $(m,n,2)$-scheme. If $m-1$ is a prime power, then $n=1$.
\elem

\prf Let $p$ be prime such that $m-1=p^t$ for some $t$. In this case, $p$ does not divide $m$ and we are done by Theorem \ref{THM:main3}~(ii).
\eprf

%\lem\label{Lemma:(m,n,m)} If $\C$ is a reduced $(m,n,m)$-scheme, $n=1$.
%\elem
%%\prf The proof is an immediate consequence of Theorem \ref{THM:main3}~(ii).
%\eprf

\lem\label{Lemma:m=2r} Let $\C$ be a reduced $(m,n,r)$-scheme and $X,Y\in\Fib(\C)$ with $X\neq Y$. If $m=2r$, then the following hold:

\begin{enumerate}
\item  [\rm~(i)]  For each $T\in\R_{X,Y}$, $d_T=2$.
\item  [\rm~(ii)] For each $R\in\R_{X}$, $d_R\in\{1,2,4\}$ and
                 $$|\{R\in\R_{X}\mid d_R=1\}|=2|\{R\in\R_{X}\mid d_R=4\}|.$$
\end{enumerate}

\elem

\prf (i) Let $\C$ be a reduced $(m,n,r)$-scheme and $X,Y\in\Fib(\C)$ with $X\neq Y$. Then as $d_T\geq2$ for each $T\in\R_{X,Y}$, it follows from
$m=\sum_{T\in\R_{X,Y}}d_T$ that $2r\leq m$ and the equality holds if and only if $d_T=2$ for each $T\in\R_{X,Y}$.\\
~(ii) Let $R\in\R_X$ and $T\in\R_{X,Y}$. Then by Lemma \ref{Lemma:constants}~(\ref{qtn:dRdS}),~(\ref{qtn:dgree-constant}), there exist none-negative integers $\alpha$ and $\beta$ such that $2d_R=d_Rd_T=\alpha d_S+\beta d_{S'}=2\alpha+2\beta$ and $\alpha, \beta\leq2$.
This implies that $d_R\leq4$. By Lemma \ref{Lemma:coprime}, $d_R\in\{1,2,4\}$.
We set $k_i:=|\{R\in\R_{X}\mid  d_R=i\}|$ for $i\in\{1,2,4\}$. Since $k_1+k_2+k_4=|\R_X|=|\R_{X,Y}|$, it follows that $m=k_1+2k_2+4k_4=2(k_1+k_2+k_4)$. Therefore, $k_1=2k_4$.
\eprf
%$\mathbb{Z}_{\geq0}$
\lem \label{Lemma:m prim Hanaki} For each $(m,n,r)$-scheme, if $m$ is prime, then $r-1$ divides $m-1$.
\elem

\prf Let $X\in\Fib(\C)$ and consider the homogeneous component $(X,\R_X)$. Since $|X|=m$ is prime, by \cite[Theorem 3.3]{Han_Uno2006} $d_R=d$ for all
$R\in\R_X$ with $R\neq\Delta_X$. Then $m-1=\sum_{\substack{R\in\R_X,\\\Delta_X\neq R}}d_R=(r-1)d.$
\eprf

\lem\label{Lemma:symmetric-odd} Let $\CVR$ be an $(m,n,r)$-scheme. If $m$ is odd, then 
each non-reflexive symmetric basis relation of $\C$ has even degree.
\elem

\prf Let $S\in\R_X$ be symmetric for some $X\in\Fib(\C)$. Since $S\neq\Delta_X$, $|S|$ is even. By (\ref{qtn:degree in&out}), $|S|=d_Sm$
and thus $d_S$ is even.
\eprf

\lem[{\rm \cite[(3.2)]{Higman-SRD95}}]\label{Lemma:Higman symmetric} Let $\C$ be a reduced $(m,n,3)$-scheme.
Then $\C_X$ is symmetric for each $X\in\Fib(\C)$.
\elem

\noindent\textbf{Proof of Theorem \ref{THM:main4}:}
So far in this section we have been preparing some lemmas which will be applied to enumerate reduced $(m,n,r)$-schemes for $m$ up to $11$. %In this case, by Theorem \ref{THM:main3}~(\ref{THM:m<2r}) we only need to enumerate $(m,n,r)$-schemes for $r$ up to $5$. In the first step, we may eliminate some cases considering $\C_X$ and then $d_X$. Note that the others will be investigated in Table \ref{Tab:enume3}.
%Now let us see how we could determine $d_X$ for a given reduced $(m,n,r)$-scheme $\C$ and $X\in\Fib(\C)$. For instance, if $\C$ is a $(7,n,r)$-scheme, then by Lemma \ref{Lemma:m prim Hanaki}, $d_X\in\set{\set{1,6},\set{1,2,2,2},\set{1,3,3}}$. Moreover, if $\C$ is reduced, then by
%Lemmas \ref{Lemma:Higman symmetric} and \ref{Lemma:symmetric-odd}, $d_X\neq\set{1,3,3}$ whereas $d_X\neq\set{1,2,2,2}$ by Theorem \ref{THM:main3}~(i).  Therefore, all of them are eliminated except $(2,7)$.

The entries $(r,m)$ such that $m<2r$ are eliminated by Theorem \ref{THM:main3}~(i) whereas $(2,m)$'s are eliminated by
Lemma \ref{Lemma:r=2} except $(2,7)$ and $(2,11)$. If $\C$ is a reduced $(m,n,2)$-scheme with $m\in\set{7,11}$, then by
Remark \ref{Remark:(m,n,2)} and Lemma \ref{Lemma:(m,n,2)a,b}, $n\leq 2$. Thus we have eliminated the first row of
Table \ref{Tab:enume1}.

Applying Lemma \ref{Lemma:m=2r} for  $(r,m)=(5,10)$ we obtain that $d_{X}=\set{1,1,2,2,4}$. According to
\cite{Nomia, Han_Miya} there are no homogeneous schemes with $d_{X}=\set{1,1,2,2,4}$ where we can prove this fact by a
theoretical way.

The entries $(4,11)$ and $(5,11)$ are eliminated by Lemma \ref{Lemma:m prim Hanaki} whereas $(3,7)$ and $(3,11)$ are eliminated by
Lemmas \ref{Lemma:Higman symmetric} and \ref{Lemma:symmetric-odd}. %(in fact there exist no $(11,1,r)$-schemes for $r\in\set{4,5}$).
%In the following table the entries being eliminated are denoted by $*$.
%Let us summarize our enumeration in the following table. The $(r,m)$-entry of the following table characterize $n$ for a reduced $(m,n,r)$-scheme.

An $(r,m)$-entry of the following table is denoted by $*$ if there exists no $(m,1,r)$-scheme.

\begin{table}[H]
\centering
\begin{tabular}{|c|c|c|c|c|c|c|c|c|cc|}
\hline
 \backslashbox{$r$}{$m$}   & 4  & 5  & 6 &  7 &   8 & 9 & 10 & 11    \\\hline
                         2 & 1  & 1  & 1 & $\leq 2$ &   1 & 1 &  1 &  $\leq$ 2  \\ \hline
                         3 & 1  & 1  & 1 &  1 & $\leq 2$ & 1 &  1 &  1    \\ \hline
                         4 & 1  & *  & 1 &  1 & $\leq 2$ & 1 &  1 &  *    \\ \hline
                         5 & *  & 1  & * &  * &   1 & 1 &  * &  *    \\
\hline
\end{tabular}
\renewcommand{\captionlabeldelim}{.}
\caption{}\label{Tab:enume1}
\end{table}
%\begin{table}[H]
%\centering
%\begin{tabular}{|c|c|c|c|c|c|c|c|c|c|}
%\hline
%  % after \\: \hline or \cline{col1-col2} \cline{col3-col4} ...
%\backslashbox{$r$}{$m$}& 4  & 5  & 6 & 7 & 8 & 9 & 10 & 11  \\ \hline
%                     2 & *  & *  & * & *  & * & * & * & *   \\ \hline
%                     3 & *  & *  &   & *  &   &   &    & *   \\ \hline
%                     4 & *  & *  & * & *  &   &   &    & *   \\ \hline
%                     5 & *  & *  & * & *  & * & * & *  & *     \\
%\hline
%\end{tabular}
%\renewcommand{\captionlabeldelim}{.}
%\caption{}\label{Tab:enume1}
%\end{table}

%The following cases are eliminated by Theorem \ref{THM:main3}~(\ref{THM:p-val}).
%\begin{table}[H]
%\centering
%\begin{tabular}{|c|c|}
%\hline
%&$(r,m)=(r,\sum_{R\in\R_X}d_R)$  \\\hline
% r=3&(9,1+4+4)\\\hline
% r=4&(8,1+1+3+3),~(9,1+2+2+4)\\ \hline
%\end{tabular}
%\renewcommand{\captionlabeldelim}{.}
%\caption{}\label{Tab:enume2}
%\end{table}
\newpage

The following is the list of $(r,m)$, $\sum_{i=1}^{r} a_i$ and $\sum_{i=1}^{r} b_i$ where $m=\sum_{i=1}^{r} a_i=\sum_{i=1}^{r} b_i$, $1=a_1\leq a_2\leq\dots\leq a_r$ and $2\leq b_1\leq b_2\dots\leq b_r$ such that $d_X=\set{a_1,\dots,a_r}$ and $d_{X,Y}=\set{b_1,\dots,b_r}$ for some $(m,1,r)$-scheme $(X,\R_X)$ not satisfying the assumption of Theorem \ref{THM:main3} (see \cite{Han_Miya, Nomia}).

\begin{table}[H]
\centering
\begin{tabular}{|l|l|c|c|}
     \hline
     $(r,m)$            & $\sum_{R\in\R_X}d_R$ & $\sum_{R\in\R_{X,Y}}d_R$& \\
     \hline

     (3,6)           &1+1+4 & 2+2+2&    Not occur by Lemma \ref{Lemma:(6,n,3)-scheme}   \\
     \hline

     (3,8)           &1+1+6& 2+2+4&    Not occur by Lemma \ref{Lemma:reduced(8,n,3)-scheme1}\\\cline{3-4}
                     &     & 2+3+3&    $n\leq2$ by Lemma \ref{Lemma:reduced(8,n,3)-scheme2} \\\cline{2-2}\cline{3-4}
                     &1+3+4& 2+2+4&    Not occur by Lemma \ref{Lemma:constants}~(\ref{qtn:symmetric-dR=2}) \\\cline{3-4}\cline{3-4}
                     &     & 2+3+3&    Not occur by Lemma \ref{Lemma:constants}~(\ref{qtn:symmetric-dR=2}) \\
     \hline

     (3,9)           &1+2+6&2+2+5& Not occur by Lemma \ref{Lemma:coprime}                    \\\cline{3-4}
                     &     &2+3+4& Not occur by Lemma \ref{Lemma:$(9,n,3)$-scheme}           \\\cline{3-4}
                     &     &3+3+3& Not occur (see Lemma \ref{Lemma:$(9,n,3)$-scheme})         \\\cline{3-4}\cline{2-2}
     \hline

     (3,10)           &1+1+8& 2+3+5&  Not occur by Lemma  \ref{Lemma:coprime}               \\\cline{3-4}
                      &     & 3+3+4&  Not occur by Lemma  \ref{Lemma:coprime}               \\\cline{3-4}
                      &     & 2+4+4&  Not occur by Lemma \ref{Lemma:$(10,n,3)$-scheme}~(i)  \\\cline{3-4}
                      &     & 2+2+6&  Not occur by Lemma \ref{Lemma:$(10,n,3)$-scheme}~(i)  \\\cline{3-4} \cline{2-2}

                      &1+3+6&2+3+5& Not occur by Lemma \ref{Lemma:coprime}                  \\\cline{3-4}
                      &     &3+3+4& Not occur by Lemma \ref{Lemma:$(10,n,3)$-scheme}~(ii)                \\\cline{3-4}
                      &     &2+4+4& Not occur by Lemma \ref{Lemma:coprime}                               \\\cline{3-4}
                      &     &2+2+6& Not occur by Lemma \ref{Lemma:constants}~(\ref{qtn:symmetric-dR=2})  \\\cline{2-2}\cline{3-4}

                      &1+4+5&2+3+5& Not occur by Lemma \ref{Lemma:coprime}   \\\cline{3-4}
                      &     &3+3+4& Not occur by Lemma \ref{Lemma:coprime}   \\\cline{3-4}
                      &     &2+4+4& Not occur by Lemma \ref{Lemma:constants}~(\ref{qtn:symmetric-dR=2})              \\\cline{3-4}
                      &     &2+2+6& Not occur by Lemma \ref{Lemma:constants}~(\ref{qtn:symmetric-dR=2})              \\\cline{3-4}
     \hline

     (4,8)            &1+1+2+4&2+2+2+2& $n\leq2$ (see Lemma \ref{Lemma:reduced (8,n,4)-scheme2})   \\
      \hline
     (4,9)            &1+1+1+6&2+2+2+3& Not occur by Lemma \ref{Lemma:constants}~(\ref{qtn:symmetric-dR=2}) and Lemma \ref{Lemma:symmetric-odd}\\\cline{2-2}\cline{3-4}
                      &1+2+3+3&2+2+2+3& Not occur by Lemma \ref{Lemma:$(9,n,4)$-scheme}                 \\
     \hline

     (4,10)          &1+2+2+5&2+2+2+4& Not occur by Lemma \ref{Lemma:coprime}   \\\cline{3-4}
                     &       &2+2+3+3& Not occur by Lemma \ref{Lemma:coprime}   \\\cline{2-2}\cline{3-4}
                     &1+1+4+4&2+2+3+3& Not occur by Lemma \ref{Lemma:coprime}   \\\cline{3-4}
                     &       &2+2+2+4& Not occur by Lemma \ref{Lemma:$(10,n,4)$-scheme} \\
    \hline
\end{tabular}
\renewcommand{\captionlabeldelim}{.}
\caption{}\label{Tab:enume3}
\end{table}

\newpage

\lem \label{Lemma:(6,n,3)-scheme} If $\C$ is a reduced $(6,n,3)$-scheme such that $d_{X}=\set{1,1,4}$ for some $X\in\Fib(\C)$, then
$d_{X,Y}\neq\set{2,2,2}$ for each $Y\in\Fib(\C)$ with $Y\neq X$.
\elem

\prf
Suppose by the contrary that $d_{X,Y}=\set{2,2,2}$ for some $Y\neq X$. By Lemma \ref{Lemma:m=2r}, $d_{Y}=\set{1,1,4}$.
Taking $R,S\in\R_{X,Y}$ with $R\neq S$ we obtain from Lemma \ref{Lemma:constants} (\ref{qtn:symmetric-dR=2}) that $R^tR=S^tS=\set{\Delta_Y,T}$ where
$T\in \R_Y$ with $T\neq\Delta_Y$ and $d_T=1$. By Lemma
\ref{Lemma:constants}~(\ref{qtn:dRdS}),~(\ref{qtn:lcm}),~(\ref{qtn:Kronecker}), $4=d_Rd_{S^t}=\alpha+4\beta$ for some non-negative integers $\alpha,\beta\leq2$. This implies $\alpha=0$ and $\beta=1$, which contradicts Lemma \ref{Lemma:constant2}.

\eprf

\lem \label{Lemma:reduced(8,n,3)-scheme1}
If $\C$ is a reduced $(8,n,3)$-scheme such that $d_X=\set{1,1,6}$ for some $X\in\Fib(\C)$, then we have the following:

\nmu[(i)]
\item For each $Y\in\Fib(\C)$ with $Y\neq X$, $d_{X,Y}\neq\set{2,2,4}$. Indeed, $d_{X,Y}=\set{2,3,3}$ for each $Y\in\Fib(\C)$.
\item Let $\R_{X,Y}=\set{R,S,S'}$ such that $d_{R}=2$ and $d_{S}=d_{S'}=3$. Let  $T\in\R_X$ with $T\neq\Delta_X$ and $d_{T}=1$.
Then $TR=\set{R}$, $TS=\set{S'}$ and $TS'=\set{S}$.
\enmu
\elem

\prf (i) Suppose by the contrary that $d_{X,Y}=\set{2,2,4}$ for some $Y\in\Fib(\C)$, and take $R\in\R_{X}$ and $S\in\R_{X,Y}$ so that
$d_R=6$ and $d_S=2$. It follows from Lemma \ref{Lemma:constants}~(\ref{qtn:dRdS}),~(\ref{qtn:lcm}),~(\ref{qtn:dgree-constant}) that for some non-negative
integers $\alpha,\beta, \gamma$ we have
$$12=d_Rd_S=2\alpha +2\beta + 4\gamma ,\quad ~6\mid2\alpha,~6\mid2\beta,~6\mid4\gamma \quad\alpha,\beta,\gamma\leq 2.$$
This implies that $\gamma=0$ and $12=2\alpha +2\beta\leq8$, a contradiction.\\
~(ii) As $d_T=1$, Lemma \ref{Lemma:constants}~(\ref{qtn:d_RS}), (\ref{qtn:gcd}) asserts that $d_{TR}=2$ and $TR=\set{R}$,
since $R$ is the unique basis relation in $\R_{X,Y}$ of degree 2.
By the same observation $d_{TS}=3$ and $TS\in\R_{X,Y}$. If $TS=\set{S}$, then by Lemma \ref{Lemma:constants}~(\ref{qtn:dRdS}), $c_{TS}^{S}=1$ and Lemma \ref{Lemma:constants}~(\ref{qtn:lcm}) implies that
$c_{SS^t}^{T}=3$. Therefore, applying Lemma \ref{Lemma:constants}~(\ref{qtn:dRdS}),~(\ref{qtn:Kronecker}) we have $9=d_Sd_{S^t}=3+3+c_{SS^t}^{T'}6$ where $T'\in\R_X$ with $d_{T'}=6$, a contradiction.
\eprf

\lem \label{Lemma:reduced(8,n,3)-scheme2}
If $\C$ is a reduced $(8,n,3)$-scheme such that $d_X=\set{1,1,6}$ for some $X\in\Fib(\C)$, then $n\leq2$.
\elem

\prf Suppose by the contrary that $X,Y$ and $Z$ are distinct fibers of $\C$. Then by Lemma \ref{Lemma:reduced(8,n,3)-scheme1},
$d_{X,Y}=d_{Y,Z}=d_{X,Z}=\set{2,3,3}$. Let $R\in\R_{X,Y}$ and $S\in\R_{Y,Z}$ with $d_R=2$ and $d_S=3$. It follows from Lemma \ref{Lemma:constants}~(\ref{qtn:dRdS}),~(\ref{qtn:lcm}),~(\ref{qtn:dgree-constant}), there exist non-negative integers $\alpha,\beta, \gamma$ such that 
$$6=d_Rd_S=2\alpha+3\beta+3\gamma,\quad 6\mid2\alpha,\quad\alpha\leq 2.$$
This implies $\alpha=0$ and $RS=\set{S'}$ where $S'\in\R_{Y,Z}$ with $d_{S'}=3$. Since $c_{RS}^{S'}=2$, by Lemma \ref{Lemma:constant2}, $R^tR\cap SS^t=\set{\Delta_Y,T}$ for some $T\in\R_Y$ with $d_T=1$. Thus $9=d_Sd_{S'}=3+3+6\alpha$. It follows that $3=6\alpha$ , a contradiction.

%Let $\R_{X,Y}=\set{R_1,S_1,{S'}_1}$, $\R_{Y,Z}=\set{R_2,S_2,{S'}_2}$ and $\R_{X,Z}=\set{R_3,S_3,{S'}_3}$ such that $d_{R_i}=2$ and $d_{S_i}=d_{{S'}_i}=3$ for $i=1,2,3$. By Lemma \ref{Lemma:constants}~(\ref{qtn:lcm}),~(\ref{qtn:dgree-constant}), $3\mid c_{S_1R_2}^{R_3}\leq 2$. Therefore, $c_{S_1R_2}^{R_3}=0$ and $R_3\notin S_1R_2$. By the same observation, $R_3\notin S_1S_2$ and $R_3\notin S_1S_2'$. Therefore,
%$S_1\R_{Y,Z}\subseteq\set{S_3,S_3'}$, which contradicts Lemma \ref{Lemma:constants}~(\ref{qtn:dgree-constant}).
\eprf

\lem\label{Lemma:$(9,n,3)$-scheme} Let $\C$  be a reduced $(9,n,3)$-scheme and $d_{X}=\set{1,2,6}$ for some fiber $X$.
Then $d_{X,Y}\notin\big\{\set{2,3,4},\set{3,3,3}\big\}$ for each $Y\in\Fib(\C)$.
\elem
\prf
Suppose that $d_{X,Y}=\set{2,3,4}$ and take the basis relations $R\in\R_X$ and $S\in\R_{X,Y}$ so that $d_R=6$ and $d_S=2$.
It follows from Lemma \ref{Lemma:constants}~(\ref{qtn:dRdS}),~(\ref{qtn:lcm}),~(\ref{qtn:dgree-constant}),
there exist non-negative integers $\alpha,\beta, \gamma$ such that
$$12=d_Rd_S=2\alpha +3\beta + 4\gamma,\quad~6\mid2\alpha,~6\mid3\beta,~6\mid4\gamma, \quad\alpha,\beta,\gamma\leq 2.$$
This implies that $\gamma=0$ and $12=2\alpha +3\beta\leq 10$, a contradiction.

Suppose that $d_{X,Y}=\set{3,3,3}$ for some $Y\in\Fib(\C)$ and take distinct $R,S\in \R_{X,Y}$. By
Lemma \ref{Lemma:constants}~(\ref{qtn:dRdS}),~(\ref{qtn:Kronecker}),
for some non-negative integers $\alpha,\beta$ we have $9=d_Rd_{S^t}=2\alpha+6\beta=2(\alpha+3\beta)$, a contradiction.
\eprf

\label{Lemma:constant3} Let $X,Y\in\Fib(\C)$ with $X\neq Y$ and $R,S,S'\in\R_{X,Y}$. Then $R^tR\cap S^tS'\neq\emptyset$ if and only if $RS^t\cap RS'^{t}\neq\emptyset$. We use this fact in the proof of the following lemma.

\lem\label{Lemma:$(10,n,3)$-scheme} Let $\C$ be a reduced $(10,n,3)$-scheme. Then the following hold:
\nmu[(i)]
\item If $d_{X}=\set{1,1,8}$ for some $X\in\Fib(\C)$, then $d_{X,Y}\notin\set{\set{2,2,6},\set{2,4,4}}$ for each $Y\in\Fib(\C)$.
\item If $d_{X}=\set{1,3,6}$ for some $X\in\Fib(\C)$, then $d_{X,Y}\neq\set{3,3,4}$ for each $Y\in\Fib(\C)$.
\enmu
\elem

\prf (i) Suppose that $d_{X,Y}=\set{2,2,6}$ for some $Y\in\Fib(\C)$. Take $R,S\in\R_{X,Y}$ with $R\neq S$ and $d_R=d_S=2$. By Lemma
\ref{Lemma:constants}~(\ref{qtn:dRdS}),~(\ref{qtn:Kronecker}),~(\ref{qtn:dgree-constant}), $4=d_Rd_{S^t}=\alpha+8\beta$ for some non-negative integers $\alpha,\beta$ with $\alpha,\beta\leq 2$. It follows that $\alpha=0$ and $4=8\beta$, a contradiction.

Suppose that $\R_{X,Y}=\set{R,S,S'}$
such that $d_R=2$ and $d_S=d_{S'}=4$. By Lemma~(\ref{qtn:dRdS}),~(\ref{qtn:Kronecker}), $8=d_Rd_{S^t}=\alpha+8\beta$ for some non-negative integers
$\alpha,\beta$ with $4\mid \alpha \leq2$. This implies that $\alpha=0$ and then $\beta=1$. Therefore, $RS^t=\set{T'}$ where $T'\in\R_X$ with $d_{T'}=8$.
By the same observation, $RS'^t=\set{T'}$. Therefore, $T\in S^tS'\cap R^tR$ where $T\neq\Delta_Y$. On the other hand,
by Lemma \ref{Lemma:constants}~(\ref{qtn:symmetric-dR=2}), $d_T=1$. It follows from Lemma \ref{Lemma:constants}~(\ref{qtn:dRdS}),~(\ref{qtn:lcm}),~(\ref{qtn:Kronecker}),~(\ref{qtn:dgree-constant}) that
for some non-negative integers $\alpha,\beta$ we have $16=d_{S^t}d_{S'}=\alpha+8\beta$ with $4\mid\alpha$ and $0<\alpha\leq4$.
This implies that $\alpha=4$ and thus $12=8\beta$,  a contradiction.

(ii) Take $R\in\R_X$ and $S\in\R_{X,Y}$ with $d_R=3$ and $d_S=4$.
By Lemma \ref{Lemma:constants}~(\ref{qtn:dRdS}),~(\ref{qtn:lcm}),~(\ref{qtn:dgree-constant}),
for some non-negative integers $\alpha,\beta,\gamma$ we have $12=d_{R}d_{S}=3\alpha+3\beta+4\gamma$ with $12\mid\alpha$ and $12\mid\beta$ and $\alpha,\beta,\gamma\leq3$. This implies that $\alpha=\beta=0$ and $\gamma=3$.
Hence $RS=\set{S}$. By Lemma \ref{Lemma:d_L_S}, $d_{L_S}\mid\rm gcd(10,4)=2$ which is a contradiction since $d_{L_S}>d_R=3$.
\eprf

\lem \label{Lemma:reduced (8,n,4)-scheme1}
Let $\C$ be a reduced $(8,n,4)$-scheme such that $d_X=\set{1,1,2,4}$ for some $X\in\Fib(\C)$. Then for all $X,Y\in\Fib(\C)$ with $X\neq Y$, there exists $R\in\R_{X,Y}$ such that $RR^t=\set{\Delta_X,S}$ \rm{(}resp. $R^tR=\set{\Delta_Y,S'}$\rm{)} where $S$ is the unique basis relation in $\R_X$ with $d_S=2$
 \rm{(}resp. $S'$ is the unique basis relation in $\R_Y$ with $d_{S'}=2$\rm{)}.
\elem

\prf Let $\C$ be a reduced $(8,n,4)$-scheme such that $d_X=\set{1,1,2,4}$ for some $X\in\Fib(\C)$. Then $d_{X,Y}=\set{2,2,2,2}$ for all
$X,Y\in\Fib(\C)$ with $X\neq Y$.
Let $T\in\R_X$ with $T\neq\Delta_X$ and $d_T=1$. Then $d_{TR}=2$ and $|TR|=1$. Suppose that $TR=\set{R}$ for each $R\in\R_{X,Y}$. Then $T\notin RS^t$ for all $R,S\in\R_{X,Y}$ with $R\neq S$. Thus by Lemma \ref{Lemma:constants}~(\ref{qtn:dRdS}),~(\ref{qtn:dgree-constant}),
$4=d_Rd_{S^t}=2\alpha+4\beta$ for some none-negative integers $\alpha,\beta$ with $\alpha\leq2$. By Lemma \ref{Lemma:constant2}, $\beta=0$ and $\alpha=2$.
This implies that $R\R_{Y,X}\subsetneq\R_X$, which contradicts Lemma \ref{Lemma:constants}~(\ref{qtn:dgree-constant}). Thus there exists $R\in\R_{X,Y}$ such that $TR\neq\set{R}$. Equivalently, $T\notin RR^t$. It follows from Lemma \ref{Lemma:constants}~(\ref{qtn:symmetric-dR=2}) that $RR^t=\set{\Delta_X, S}$ where $S$ is the unique basis relation in $\R_X$ with $d_S=2$
\eprf

\lem\label{Lemma:reduced (8,n,4)-scheme2}If $\C$ is a reduced $(8,n,4)$-scheme such that $d_X=\set{1,1,2,4}$ for some $X\in\Fib(\C)$, then $n\leq 2$.
\elem

\prf
Suppose by the contrary that $X,Y$ and $Z$ are distinct fibers of $\C$. Then by Lemma \ref{Lemma:reduced (8,n,4)-scheme1}, there exist $R\in\R_{X,Y}$ and $T\in\R_{Y,Z}$ such that $R^tR=TT^t=\set{\Delta_X, S}$ where $S$ is the unique basis relation in $\R_Y$ with $d_S=2$. It follows from
Lemma \ref{Lemma:constants}~(\ref{qtn:dRdS}) that $c_{TT^t}^S=c_{R^tR}^S=1$. Let $(y,y')\in S$. Then there exists $(x,z)\in X\times Z$ such that
$R_{in}(y)\cap R_{in}(y')=\set{x}$ and $T_{out}(y)\cap T_{out}(y')=\set{z}$. As $d_T=2$, we may assume that $T_{out}(y)=\set{z,z_1}$ and
$T_{out}(y')=\set{z,z_2}$. Note that $z_1\neq z_2$, otherwise $c_{TT^t}^S \geq2$, a contradiction. This means that $(R\circ T)_{out}(x)=\set{z,z_1,z_2}$
and thus $d_{RT}=d_{R\circ T}=3$, which is a contradiction since $d_{RT}$ must be a sum of degrees in $d_{X,Z}=\set{2,2,2,2}$.
\vspace{-3mm}
\begin{figure}[H]
\begin{displaymath}
    \xymatrix{
            &\VRT{x}\ar[dl]_{\scriptscriptstyle\mathit{R}} \ar[dr]^{\scriptscriptstyle\mathit{R}}& \\
            \VRT{y}\ar[d]_{{\scriptscriptstyle\mathit{T}}}\ar[dr]^{{\scriptscriptstyle\mathit{T}}} \ar[rr]_{{\scriptstyle\mathit{s}}}&  & \VRT{y'}\ar[dl]_{{\scriptscriptstyle\mathit{T}}}\ar[d]^{{\scriptscriptstyle\mathit{T}}} \\
            \VRT{z_1}& \VRT{z} & \VRT{z_2}}
\end{displaymath}\label{fig:6}
\vspace{-4mm}
\renewcommand{\captionlabeldelim}{.}
\caption{}\label{Fig:(8,n,4)-scheme}
\end{figure}
\eprf

\exm The association scheme \verb"as16 No.122" as in {\rm \cite{Han_Miya}} induces the thin residue fission  {\rm(see \cite[Proposition 3.1]{28points})}, which is a reduced $(8,2,4)$-scheme whose relational matrix is
$$
\left(
  \begin{array}{cccccccc|cccccccc}
0 & 1 & 2 & 2 & 3 & 3 & 3 & 3 & 4 & 4 & 5 & 5 & 6 & 7 & 6 & 7 \\
1 & 0 & 2 & 2 & 3 & 3 & 3 & 3 & 4 & 4 & 5 & 5 & 7 & 6 & 7 & 6 \\
2 & 2 & 0 & 1 & 3 & 3 & 3 & 3 & 5 & 5 & 4 & 4 & 6 & 7 & 7 & 6 \\
2 & 2 & 1 & 0 & 3 & 3 & 3 & 3 & 5 & 5 & 4 & 4 & 7 & 6 & 6 & 7 \\
3 & 3 & 3 & 3 & 0 & 1 & 2 & 2 & 6 & 7 & 6 & 7 & 4 & 4 & 5 & 5 \\
3 & 3 & 3 & 3 & 1 & 0 & 2 & 2 & 7 & 6 & 7 & 6 & 4 & 4 & 5 & 5 \\
3 & 3 & 3 & 3 & 2 & 2 & 0 & 1 & 6 & 7 & 7 & 6 & 5 & 5 & 4 & 4 \\
3 & 3 & 3 & 3 & 2 & 2 & 1 & 0 & 7 & 6 & 6 & 7 & 5 & 5 & 4 & 4 \\\hline
4' & 4' & 5' & 5' & 6' & 7' & 6' & 7' & 0' & 1' & 2' & 2' & 3' & 3' & 3' & 3' \\
4' & 4' & 5' & 5' & 7' & 6' & 7' & 6' & 1' & 0' & 2' & 2' & 3' & 3' & 3' & 3' \\
5' & 5' & 4' & 4' & 6' & 7' & 7' & 6' & 2' & 2' & 0' & 1' & 3' & 3' & 3' & 3' \\
5' & 5' & 4' & 4' & 7' & 6' & 6' & 7' & 2' & 2' & 1' & 0' & 3' & 3' & 3' & 3' \\
6' & 7' & 6' & 7' & 4' & 4' & 5' & 5' & 3' & 3' & 3' & 3' & 0' & 1' & 2' & 2' \\
7' & 6' & 7' & 6' & 4' & 4' & 5' & 5' & 3' & 3' & 3' & 3' & 1' & 0' & 2' & 2' \\
6' & 7' & 7' & 6' & 5' & 5' & 4' & 4' & 3' & 3' & 3' & 3' & 2' & 2' & 0' & 1' \\
7'' & 6' & 6' & 7'' & 5' & 5' & 4' & 4' & 3' & 3' & 3' & 3' & 2' & 2' & 1' & 0' \\
  \end{array}
\right).
$$

Also the thin residue fission of the association scheme \verb"as16 No.51"  as in {\rm\cite{Han_Miya}}, is a reduced $(8,2,3)$-scheme.
\eexm

\lem\label{Lemma:$(9,n,4)$-scheme} Let $\C$ be a reduced $(9,n,4)$-scheme such that $d_{X}=\set{1,2,3,3}$ for some $X\in\Fib(\C)$. Then $d_{X,Y}\neq\set{2,2,2,3}$ for each $Y\in\Fib(\C)$.
\elem

\prf Suppose by the contrary that $R_1,R_2,R_3\in \R_{X,Y}$ with $d_{R_i}=2$ for $i\in\set{1,2,3}$. For all $i,j\in\set{1,2,3}$ with $i\neq j$, by
Lemma \ref{Lemma:constants}~(\ref{qtn:dRdS}),~(\ref{qtn:Kronecker}) we have $4=d_{R_i}d_{R_j^t}=2\alpha+3\beta+3\gamma$. This implies that $\beta=\gamma=0$ and $\alpha=2$. Hence, for all $i,j\in\set{1,2,3}$ with $i\neq j$, $R_iR_j^t=\set{T}$ where $T\in\R_X$ with $d_T=2$. It follows from Lemma \ref{Lemma:constants}~(\ref{qtn:dT<dR}) that $\set{R_1}=TR_2=\set{R_3}$, a contradiction.
\eprf

\lem\label{Lemma:$(10,n,4)$-scheme} Let $\C$ be a reduced $(10,n,4)$-scheme such that $d_{X}=\set{1,1,4,4}$ for some $X\in\Fib(\C)$.
Then $d_{X,Y}\neq\set{2,2,2,4}$ for each $Y\in\Fib(\C)$.
\elem

\prf Suppose by the contrary that $d_{X,Y}=\set{2,2,2,4}$ for some $Y\in\Fib(\C)$ with $Y\neq X$. According to \cite{Nomia, Han_Miya},
$d_Y\in\set{\set{1,1,4,4},\set{1,2,2,5}}$. It follows from Lemma \ref{Lemma:coprime} that $d_Y=\set{1,1,4,4}$. Take $R,S\in \R_{X,Y}$ with $R\neq S$ and $d_{R}=d_S=2$. By Lemma \ref{Lemma:constants}~(\ref{qtn:dRdS}),~(\ref{qtn:Kronecker}),~(\ref{qtn:dgree-constant}), there exist non-negative
integers $\alpha,\beta,\gamma$ such that
$$4=d_Rd_{S^t}=\alpha +4\beta+4\gamma,\quad  \alpha,\beta, \gamma\leq 2.$$
This implies that $4\mid\alpha$. Hence, $\alpha=0$ and $\beta+\gamma=1$ which contradicts Lemma \ref{Lemma:constant2}.
%By Lemma \ref{Lemma:constants}~(\ref{qtn:symmetric-dR=2}), there exists $T\in R^tR\cap S^tS$
%with $T\neq\Delta_Y$. It follows from Lemma \ref{Lemma:constant2} and Lemma \ref{Lemma:constants}, $c_{RS^t}^{T'}\geq 2$ for some $T'\in\R_{X}$.
\eprf

\noindent\textbf{Acknowledgements}\\
The authors would like to thank I. Ponomarenko, E. Bannai and A. Ivanov for their valuable comments to this article, especially internal direct sums, $(m,2,r)$-schemes and tensor products, respectively.

\end{document}